\documentclass{gtart}
\usepackage{amsmath,amscd,amssymb}
\input gtmonout    
\volumenumber{1}
\volumeyear{1998}
\volumename{The Epstein birthday schrift}
\pagenumbers{341}{364}
\received{20 November 1997}
\published{30 October 1998}
\papernumber{17}

\newtheorem{theorem}{Theorem}[section]
\newtheorem{prop}[theorem]{Proposition}

\newtheorem{cor}[theorem]{Corollary}

\theoremstyle{definition}

\newtheorem*{defn}{Definition}

\begin{document}

\title{Coarse extrinsic geometry: a survey}

\author{Mahan Mitra}

\address{Institute of Mathematical Sciences, C.I.T. Campus\\
Madras (Chennai) -- 600113, India}
\email{mitra@imsc.ernet.in}

\begin{abstract}
This paper is a survey of some of the developments in coarse
extrinsic geometry since its inception in the work of Gromov.
Distortion, as measured by comparing the diameter of balls
relative to different metrics, can be regarded as one of the
simplist extrinsic notions. Results and examples concerning
distorted subgroups, especially in the context of hyperbolic
groups and  symmetric spaces, are exposed. Other topics
considered are quasiconvexity of subgroups; behaviour at
infinity, or more precisely continuous extensions of embedding maps
to Gromov boundaries in the context of hyperbolic groups  acting
by isometries on hyperbolic metric spaces; and distortion as
measured using various other filling invariants.
\end{abstract}



\maketitle

\cl{\small\it To David Epstein on his sixtieth birthday}

\section{Introduction}

Extrinsic geometry deals with the study of the geometry
of subspaces relative to that of an ambient space. Given
a Riemannian manifold $M$ and a submanifold $N$, classical
(differential) extrinsic geometry studies infinitesimal
changes in the 
Riemannian
metric on $N$  induced from $M$. This
involves an analysis of the second fundamental form
or shape operator \cite{Grsig}. In coarse geometry local
or infinitesimal machinery is absent. Thus it does not
make sense to speak of tangent spaces or Riemannian
metrics. However, the large scale notion of metric
continues to make sense. Given a metric space $X$ and 
a subspace $Y$ one can still compare the intrinsic metric on $Y$
to the metric inherited from $X$.
This is especially useful for finitely generated subgroups of
finitely generated groups.
To formalize this, Gromov introduced the notion of distortion in
his seminal paper \cite{Gromov2}.
\begin{defn}{\rm (\cite{Gromov2},\cite{farb})}\stdspace
If $i\co \Gamma_H \rightarrow \Gamma_G$ is an embedding of the Cayley
graph of $H$ into that of $G$, then the {\it distortion} function is
given by
\begin{center}
$disto(R) = Diam_{\Gamma_H}({\Gamma_H}{\cap}B(R))$,
\end{center}
where $B(R)$ is the ball of radius $R$ around $1\in{\Gamma_G}$.
\end{defn}
The definition above differs from the one in \cite{Gromov2} by
a linear factor and coincides with that in \cite{farb}.

{\bf Note}\stdspace The above definition continues to make sense when $\Gamma_G$
and $\Gamma_H$ are replaced by graphs or (more generally)
path-metric spaces (see below for definition) $X$ and $Y$ respectively. 
\begin{defn} A {\em path-metric space} is a metric space $(X,d)$
such that for all $x, y \in X$ there exists an isometric
embedding $f \co [0,d(x,y)] \rightarrow X$ with $f(0) = x$
and $f(d(x,y)) = y$. 
\end{defn}
If the distortion function is linear we say $\Gamma_H$ (or $Y$)
is {\it quasi-isometrically}  (often abbreviated to
{\it qi\/}) embedded in $\Gamma_G$ (or $X$).
This is equivalent to the following:
\begin{defn}
A map $f$ from one metric space $(Y,{d_Y})$ into another metric space 
$(Z,{d_Z})$ is said to be
 a {\it $(K,\epsilon)$--quasi-isometric embedding} if
 
\begin{center}
${\frac{1}{K}}({d_Y}({y_1},{y_2}))-\epsilon\leq{d_Z}(f({y_1}),f({y_2}))\leq{K}{d_Y}({y_1},{y_2})+\epsilon.$
\end{center}

If  $f$ is a quasi-isometric embedding, and every point of $Z$ lies at
a uniformly bounded distance from some $f(y)$ then $f$ is said to be a
{\it quasi-isometry}. A $(K,{\epsilon})$--quasi-isometric embedding
that is a quasi-isometry will be called a\break $(K,{\epsilon})$--quasi-isometry.
\end{defn}
We collect here a few other closely related notions:
\begin{defn}A subset $Z$ of $X$ is said to be {\it $k$--quasiconvex}
if any geodesic joining $a,b\in Z$ lies in a $k$--neighborhood of $Z$.
A subset $Z$ is {\it quasiconvex} if it is $k$--quasiconvex for some $k$.

A {\it $(K,\epsilon)$--quasigeodesic} is a
$(K,\epsilon)$--quasi-isometric embedding of a closed interval in
$\mathbb{R}$. A $(K,0)$--quasigeodesic will also be called a $K$--quasigeodesic.
\end{defn}

For hyperbolic metric spaces (in the sense of Gromov \cite{Gromov})
the notions of quasiconvexity and qi embeddings coincide. This
is because quasigeodesics lie close to geodesics in hyperbolic
metric spaces \cite{Shortetal}, \cite{GhH}, \cite{CDP}.

Distortion can be regarded, in some sense, as the simplest extrinsic
notion in coarse geometry. However a complete understanding of
distortion is lacking even in special situations like subgroups
of hyperbolic groups or discrete (infinite co-volume) subgroups
of higher rank semi-simple Lie groups. One of the aims of this
survey is to expose some of the issues involved. This is done in
Section~2. 

A characterisation of quasi-isometric embeddings in terms of group
theory is another topic of extrinsic geometry that has received
some attention of late. This will be dealt with in Section 3.

A different perspective of coarse extrinsic geometry comes from the
asymptotic point of view. The issue here is behavior `at infinity'.
 From this perspective it seems possible to introduce and study
finer invariants involving distortion along specified directions.
Section 4 deals with this in the special context of hyperbolic
subgroups of hyperbolic groups.

Finally in Section 5, we discuss some other invariants of extrinsic
geometry that have come up in different contexts. 

It goes without saying that this survey reflects the author's bias
and is far from comprehensive.

{\bf Acknowledgements}\qua I would like to thank the
organizers  of David Epstein's birthday fest where the idea
of this survey first came up. I would also like to take this
opportunity to thank Benson Farb for several inspiring 
discussions on coarse geometry over the years.

\section{Distortion}

If a finitely generated subgroup $H$ of a finitely generated group $G$
is qi--embedded  we shall refer to it as undistorted.
Otherwise $H$ will be said to be distorted. We shall also have
occasion
to replace the Cayley graph of $G$ by a symmetric space (equipped with
its invariant metric) or more generally a path metric space $(X,d)$. In the
latter case, distortion will be measured with respect to the metric
$d$ on $X$.

Distorted subgroups of hyperbolic groups or symmetric spaces are
somewhat difficult to come by. This has resulted in a limited
supply of examples. Brief accounts will be given of some of the
known sources of examples. 

An aspect that will not be treated in any detail is the connection
to algorithmic problems, especially the Magnus problem. See
 \cite{Gromov2} or (for a more detailed account)
 \cite{farb} for a treatment.

\sh{Subgroups of hyperbolic groups and $SL_2{({\mathbb{C}})}$}

One of the earliest classes of examples of distorted hyperbolic
subgroups
of hyperbolic groups came from Thurston's work on $3$--manifolds 
fibering over the circle \cite{Thurston}. Let $M$ be a closed
hyperbolic $3$--manifold fibering over the circle with fiber $F$. Then 
${\pi_1}(F)$ is a hyperbolic subgroup of the hyperbolic
group ${\pi_1}(M)$. The distortion is easily seen to be 
exponential. 

It follows from work of Bonahon \cite{Bonahon1} and
Thurston \cite{Thurstonnotes} that if $H$ is a closed
surface subgroup of the fundamental group ${\pi_1}(M)$
of a closed hyperbolic $3$--manifold $M$ then the 
distortion of $H$ is either linear or exponential.
This continues to be true if $H$ is replaced by
any freely indecomposable group. In fact exponential
distortion of a freely indecomposable group corresponds
precisely (up to passing to a finite cover of $M$)
to the case of a hyperbolic $3$--manifold fibering over the
circle.

The situation is considerably less clear when we come to freely
decomposable subgroups of hyperbolic $3$--manifolds. The tameness
conjecture (attributed to Marden \cite{marden}) asserts that the covering of
a closed hyperbolic $3$--manifold corresponding to a finitely
generated subgroup of its fundamental group is topologically
tame, ie is homeomorphic to the interior of a compact $3$--manifold with boundary. If this conjecture were true, it
would follow (using a Theorem of Canary \cite{canary})
that any finitely generated subgroup $H$
of the fundamental group ${\pi_1}(M) = G$ is either
quasiconvex in $G$ or is exponentially distorted.
Moreover, exponential
distortion corresponds
precisely (up to passing to a finite cover of $M$)
to the case of a fiber of a
hyperbolic $3$--manifold fibering over the
circle. Much of this theory can be extended to take parabolics
into account.

This class of examples can be generalized in two directions. One can
ask for distorted discrete subgroups of $SL_2{(\mathbb{C})}$
or for distorted hyperbolic subgroups of hyperbolic groups (in
the sense of Gromov).  We look first at discrete subgroups
of $SL_2{(\mathbb{C})}$. A substantial class of examples comes from
geometrically tame groups. In fact the simplest surface group,
the fundamental group of a punctured torus (the puncture corresponds
to a parabolic element), displays much of the exotic extrinsic
geometry
that may occur. These examples were studied in great detail by Minsky
in
\cite{Minsky1}. The distortion function was calculated in \cite{mitra3}.

Let $S$ be a hyperbolic punctured torus so that the two shortest geodesics
$a$ and $b$ are orthogonal and of equal length. Let $S_0$ denote $S$ minus
a neighborhood of the cusp. Let $N_{\delta}(a)$ and
$N_{\delta}(b)$ be regular collar neighborhoods of $a$ and $b$ in $S_0$.
For $n\in{\mathbb{N}}$, define ${\gamma_n} = a$ if $n$ is even and equal
to $b$ if $n$ is odd. Let $T_n$ be the open solid torus neighborhood
of ${\gamma_n}{\times}\{n+{\frac{1}{2}}\}$ in ${S_0}{\times}[0,{\infty})$
given by 
\begin{center}
$T_n = N_{\delta}({\gamma_n}){\times}(n,n+1)$
\end{center}
and let
$M_0 = ({S_0}){\times}[0,{\infty}){\setminus}{\bigcup}_{n\in{\mathbb{N}}}T_n$.

Let $a(n)$ be a sequence of positive integers greater than one. Let
${\hat{\gamma_n}} = {\gamma_n}{\times}\{{n}\}$ and let $\mu_n$ be an
oriented  meridian for $\partial{T_n}$ with a single positive intersection
with ${\hat{\gamma_n}}$. Let $M$ denote the result of gluing to each 
$\partial{T_n}$  a solid torus $\hat{T_n}$,
such that the curve ${{\hat{\gamma_n}}^{a(n)}}{\mu_n}$ is glued to a 
meridian.  Let
$q_{nm}$ be the mapping class 
from $S_0$ to itself obtained by identifying $S_0$ to ${S_0}{\times}m$,
pushing through $M$ to ${S_0}{\times}n$ and back to $S_0$. Then
$q_{n(n+1)}$ is given by $\Phi_n = 
D_{\gamma_n}^{a(n)}$, where $D_c^k$ denotes Dehn 
twist along $c$, $k$ times. Matrix representations of $\Phi_n$ are
given by 
\[ \Phi_{2n} = \left( \begin{array}{cc}
                         1 & a(2n) \\ 0 & 1
                   \end{array} \right)  \]
and \[ \Phi_{2n+1} = \left( \begin{array}{cc}
                         1 & 0 \\ a(2n+1) & 1
                   \end{array} \right).  \]

Recall that the metric on $M_0$ is the restriction of the product metric.
The $\hat{T_n}$'s are given hyperbolic metrics such that their boundaries
are uniformly quasi-isometric to $\partial{T_n} \subset M_0$. Then from
\cite{Minsky1}, $M$ is quasi-isometric to the complement of a rank
one cusp in the convex core of 
a hyperbolic manifold $M_1 = {{\mathbb{H}}^3}/{\Gamma}$.
 Let $\sigma_n$ denote the shortest path
from $S_0{\times}1$ 
to $S_0{\times}n$. Let ${\overline{\sigma_n}}$ denote $\sigma_n$ with reversed
orientation. Then $\tau_n = \sigma_n{\gamma_n}{\overline{\sigma_n}}$
is a closed path in $M$ of length $2n+1$. 
Further $\tau_n$ is homotopic to a curve 
$\rho_n = \Phi_1{\cdots}\Phi_n{(\gamma_n)}$ on $S_0$. Then
\begin{center}
$\Pi_{i=1{\cdots}n}{a(i)} \leq l({\rho_n}) \leq \Pi_{i=1{\cdots}n}{(a(i)+2)}$
\end{center}

Hence
\begin{center}
$\Pi_{i=1{\cdots}n}{a(i)} \leq (2n+1)disto(2n+1) \leq
  \Pi_{i=1{\cdots}n}{(a(i)+2)}$
\end{center}

Since $M$ is quasi-isometric to the complement of the cusp of a hyperbolic
manifold \cite{Minsky1} 
and $\gamma_n$'s lie in a complement of the cusp, the distortion
function of $\Gamma$ is of the same order as the distortion function above.
In particular, functions of arbitrarily fast growth may be
realized. This
answers a question posed by Gromov in \cite{Gromov2} page 66.

A closely related class of examples (the so called `drill--holes'
examples of which the punctured torus examples above 
may be regarded as  special
cases) appears in work of Thurston \cite{Thurston} and
Bonahon and Otal \cite{bon-otal}.

Let us now turn to finitely generated subgroups of hyperbolic
groups. If we restrict ourselves to hyperbolic subgroups
there is a considerable paucity of examples. The chief
ingredient for constructing distorted hyperbolic subgroups
of hyperbolic groups is the celebrated combination theorem
of Bestvina and Feighn \cite{BF}. This theorem was partly
motivated by Thurston's hyperbolization theorem for Haken
manifolds \cite{ctm}, \cite{Thurston} and continues to be
an inevitable first step in constructing any distorted 
hyperbolic subgroups.  The following Proposition summarizes
these examples. The proof follows easily from normal
forms.

\begin{prop}
 Let $G$ be a hyperbolic group acting cocompactly on a simplicial tree
$T$ such that all vertex and edge stabilizers are hyperbolic. Also
suppose that every inclusion of an edge stabilizer in a vertex stabilizer
is a quasi-isometric embedding. Let $H$ be the stabilizer of a vertex or
edge of $T$. Then the distortion of $H$ is linear or exponential.
\label{BF}
\end{prop}

Based on Bestvina and Feighn's combination theorem and work
of Thurston's on stable and unstable foliations of surfaces \cite{FLP},
Mosher \cite{Mosher2} constructed a class of examples of normal
surface subgroups of hyperbolic groups where the quotient is
free of rank strictly greater than one.

This idea was used by Bestvina, Feighn and Handel in \cite{BFH} 
to construct similar examples where the normal subgroup is free.

Thus one has examples of exact sequences

\begin{center}

$1 \rightarrow N \rightarrow G \rightarrow Q \rightarrow 1$

\end{center}

of hyperbolic groups where $N$ is a free group or a surface group.
Owing to a general theorem of Mosher's regarding the existence
of quasi-isometric sections of $Q$ \cite{Mosher} the distortion
of any normal hyperbolic subgroup $N$ of infinite index in
a hyperbolic group $G$
is exponential.

Further, it follows from
work of Rips and Sela \cite{Sela}, \cite{rips-sela} that a torsion
free normal hyperbolic subgroup of a hyperbolic group is a free
product of free groups and surface groups. However, the only known
restriction on $Q$ is that it is hyperbolic \cite{Mosher}. It seems
natural to wonder if
there exist  examples where the exact sequence does not split or at 
least where $Q$ is not virtually free.

We now describe some examples exhibiting higher distortion \cite{mitra3}.
Start with  a hyperbolic group $G$ such that $1 \rightarrow F \rightarrow G
\rightarrow F \rightarrow 1$ is exact, where $F$ is free of rank 3.

Let $F_1 \subset G$ denote the normal subgroup. Let $F_2 \subset G$ denote
a section of the quotient group. Let $G_1,{\cdots},G_n$ be $n$ distinct
copies of $G$. Let $F_{i1}$ and $F_{i2}$ denote copies of $F_1$ and $F_2$
respectively in $G_i$. Let \\ \begin{center}
$G = {G_1}{*_{H_1}}{G_2}*{\cdots}{*_{H_{n-1}}}{G_n}$
\end{center}
where each $H_i$ is a free group of rank 3, the image of $H_i$ in $G_i$
is $F_{i2}$ and the image of $H_i$ in $G_{i+1}$ is $F_{(i+1)1}$. Then $G$
is hyperbolic. 

Let $H = F_{11} \subset G$. Then the distortion of $H$ is superexponential
for $n > 1$. In fact, it can be  checked inductively 
that the distortion function 
is an iterated exponential of height $n$.

Starting from Bestvina, Feighn and Handel's examples above, one
can construct examples with distortion a tower function.
Let ${a_1},{a_2},a_3$ be 
generators of $F_1$ and ${b_1},{b_2},{b_3}$ be generators of $F_2$. Then
\begin{center}
$G = \{ {a_1},{a_2},{a_3},{b_1},{b_2},{b_3} : {b_i^{-1}}{a_j}{b_i} = w_{ij} \}$
\end{center}
where $w_{ij}$ are words in $a_i$'s. We add a letter $c$ conjugating $a_i$'s
to `sufficiently random' words in $b_j$'s to get $G_1$. Thus,
\begin{center}
$G_1 = \{ {a_1},{a_2},{a_3},{b_1},{b_2},{b_3},c : {b_i^{-1}}{a_j}{b_i} = w_{ij}, {c^{-1}}{a_i}c = v_i \}$,
\end{center}
where $v_i$'s are words in $b_j$'s satisfying a small-cancellation type condition to ensure that $G_1$ is hyperbolic. See \cite{Gromov}, page 151 for details
on addition of `random' relations.

It can be checked that these examples have distortion function greater than
any iterated exponential.

The above set of examples were motivated largely by examples of distorted
cyclic subgroups in \cite{Gromov2}, page 67 and \cite{gervol1} (these
examples will be discussed later in this paper). 

So far, there is no satisfactory way of manufacturing examples of hyperbolic
subgroups of hyperbolic groups exhibiting arbitrarily high distortion.
It is easy to see that a subgroup of sub-exponential distortion is
quasiconvex \cite{Gromov2}. Not much else is known. 
One is thus led to the following question:
\par\medskip
{\bf Question}\stdspace Given any increasing function $f \co \mathbb{N}
\rightarrow \mathbb{N}$ does there exist a hyperbolic subgroup
$H$ of a hyperbolic group $G$ such that the distortion of
$H$ is of the order of $e^{f(n)}$?
\par\medskip
Note that the above question has a positive answer if $G$ is replaced
by $SL_2({\mathbb{C}})$.

If one does not restrict oneself to hyperbolic subgroups of hyperbolic
groups, one has a large source of examples coming from 
finitely generated subgroups of small cancellation groups.
These examples are due to Rips \cite{rips}.

Let $Q = \{ {g_1}, {\cdots}, {g_n} : {r_1},{\cdots},{r_m} \}$
be any finitely presented group.
Construct a small cancellation ($C^{\prime}(1/6)$) group
$G$ with presentation as follows:

\begin{center}

$G = \{ {g_1},{\cdots},{g_n},{a_1},{a_2} : {g_i^{-1}}{a_j}{g_i} =
u_{ij}, {g_i}{a_j}{g_i^{-1}} = v_{ij}, r_k = w_k $ \\
for $i = 1 \cdots n$, $j = 1, 2$ and $k = 1 \cdots m$. $ \} $

\end{center}

where $u_{ij}$, $v_{ij}$, $w_k$ are words in $a_1$, $a_2$
satisfying $C^{\prime}(1/6)$.

Then one has an exact sequence $ 1 \rightarrow H \rightarrow G
\rightarrow Q \rightarrow 1$ where $H$ is the subgroup of $G$
generated by ${a_1},{a_2}$ and $Q$ is the given finitely presented
group. The distortion of $H$ can be made to vary by varying $Q$
(one basically needs to vary the complexity of the word problem in
$Q$). However the subgroups $H$ are generally not finitely presented.

A remarkable example of a finitely presented normal subgroup $H$
of a hyperbolic
group $G$ has recently been discovered by Brady \cite{brady}. This 
is the first example of a finitely presented non-hyperbolic
subgroup of   a hyperbolic group. The distortion in this example
is exponential as the quotient group is infinite cyclic.

\sh{Distortion in symmetric spaces}

Now let $G$ be a semi-simple Lie group. Cyclic discrete subgroups
generated by unipotent elements are exponentially distorted. This
is because discrete subgroups of the nilpotent subgroup $N$
in a $KAN$ decomposition of
$G$ is distorted in this way. 
 This is the most well known source of distortion. 

Other known  examples seem to have their origin in rank 1 
phenomena. Given any Lie group $G$ containing ${F_2}{\times}{F_2}$
as a discrete subgroup one has distorted subgroups coming
from a construction due to Mihailova \cite{mihailova}, \cite{Gromov2},
\cite{farb} (see below). In some sense these examples are `reducible'.
Truly higher rank phenomena are hard to come by. One has the
following basic question:
\par\medskip
{\bf Question}\stdspace Are there examples of distorted finitely generated
discrete subgroups $H$ of irreducible lattices in higher rank semi-simple
Lie groups $G$ such that $H$ has no unipotent element?
(See \cite{farb} also).
\par\medskip
Note that Thurston's construction of normal subgroups cannot
possibly go through here on account of the following basic
theorem of Kazhdan--Margulis:

\begin{theorem}{\rm\cite{margulis}}\qua
Let $\Gamma$ be an irreducible lattice in a symmetric space
of real rank greater than one. Then any normal subgroup $\Lambda$
of $\Gamma$ is either finite or the quotient $\Gamma{/}\Lambda$
is finite.
\end{theorem}

Another non-distortion theorem has recently been proven by
Lubotzky--Moses--Raghunathan \cite{LMR}
answering a question of Kazhdan:

\begin{theorem} 
Any irreducible lattice in a symmetric space $X$ of rank greater than one 
is undistorted in $X$. 
\end{theorem}

The above theorems indicate the difficulty in obtaining 
distorted subgroups of higher rank Lie groups.

Similar questions may be asked for rank one symmetric spaces also
eg for complex hyperbolic, quaternionic hyperbolic
and the Cayley hyperbolic planes. Here, too there is a dearth
of examples.

In real hyperbolic spaces, the situation is slightly better
owing to Thurston's examples of $3$--manifolds fibering over the
circle. Based partly on Thurston's examples, Bowditch and
Mess \cite{bow-mess}  have described
an example of a finitely generated subgroup of a uniform lattice
in $SO(4,1)$ that is not finitely presented. 
 Abresch and Schroeder \cite{ab-sch} have given an arithmetic
construction of this lattice, too. One wonders
if this arithmetic description can be used to give similar
examples in $SU(4,1)$ or $Sp(4,1)$. 

Such infinitely presented  subgroups are
 necessarily distorted. Related examples have also been
discovered by Potyagailo and Kapovich \cite{potyagailo},
\cite{pot-kap}.

A natural question is whether Thurston's construction goes through
in higher dimensions or not:
\par\medskip
{\bf Question}\stdspace Does there exist a uniform lattice in a rank one
symmetric (other than ${\mathbb{H}}^3$)
space containing a finitely presented
(or even finitely generated) infinite
normal subgroup of infinite index?
\par\medskip
One should note that any such normal subgroup 
cannot be hyperbolic (by \cite{rips-sela}).

\sh{Distortion in finitely presented groups}

There are certain special classes of distorted subgroups of finitely
presented groups that do not fall into any of the above categories.

A basic class of examples comes from the Baumslag Solitar groups

\begin{center}
$BS(1,n) = \{ a, t : tat^{-1} = a^n \}$
\end{center}

where the cyclic group generated by $a$ has exponential distortion
for $n > 1$.

A class of examples with higher distortion have appeared in work
of Gersten \cite{gervol1}. We briefly describe these.

Take $G = \{ {g_1},{\cdots},{g_n} : {g_{i-1}}^{g_i} = {g_{i-1}^2}$
for $i = 2 \cdots n \}$.  Then the cyclic subgroup generated by
$g_1$ has distortion an iterated exponential function of height 
$n$.

Next consider $G = \{ a, b, c: a^b = a^2, a^c = b \}$. Then the 
cyclic group generated by $a$ has distortion greater than any iterated
exponential. 

Another class of subgroups with distortion a fractional power occurs
in work of Bridson \cite{bridson}: \\
Let ${G_c} = {\mathbb{Z}}^c {\rtimes}_{\phi_c} \mathbb{Z}$
where $\phi_c \in GL_c{\mathbb{Z}}$ is the unipotent matrix with
ones on the diagonal and superdiagonal and zeroes elsewhere. 
For $c > 1$, $G_c$ has infinite cyclic center. Given two such
groups  $G_a$, $G_b$ amalgamate them along their cyclic center
$\langle z\rangle $ to get $G(a,b) = {G_a}{*_{\langle z\rangle }}{G_b}$. Then the distortion
function of $G_b$ in $G(a,b)$ is of the form $n^{\frac{a}{b}}$.

A large class of examples of distortion arise from subgroups of
nilpotent and solvable groups \cite{Gromov2}.

Finally we describe a class of examples due to Mihailova \cite{mihailova}
which give rise to non-recursive distortion (see also \cite{Gromov2}
\cite{farb}). Let $G = \{ {g_1},{\cdots},{g_n} : {r_1}{\cdots}{r_m}\}$
 be any finitely presented group with
defining presentation $f \co F_n \rightarrow G$. Then $f \times f$
maps $F_n \times  F_n$ to $G \times G$. The pull-back $H$
under this
map of the `diagonal subgroup' $\{ (g,g) : g \in G \}$ is generated
by elements of the form $({g_i},{g_i})$, $i = 1 \cdots n$  and  
$(1,{r_j})$, $j = 1 \cdots m$. If $G$ has unsolvable word problem,
then the distortion of $H$ in ${F_n} \times F_n$ is non-recursive.

\section{Characterization of quasiconvexity}

It was seen in the previous section that construction of distorted
subgroups usually involves some amount of work. In fact for subgroups
of hyperbolic groups, Gromov \cite{Gromov} describes `length--angle'
relationships between generators that would ensure quasiconvexity
of the subgroup.  This can be taken as a genericity result. In another
setting, one could ask for examples of groups all whose finitely generated
subgroups are undistorted. This is known for free groups, surface groups
and abelian groups. 

However, a general group-theoretic characterization of quasiconvexity
seems far off.
Gersten has recently described a functional analytic approach to this problem.
We briefly describe this. Later we shall discuss a more group-theoretic
approach. We shall restrict ourselves to finitely generated subgroups
of hyperbolic groups (in the sense of Gromov)  in this section.

The following discussion appears in \cite{gersten1}, \cite{gersten2},
\cite{gersten3}. Let $X^{\prime}$ be a complex of type $K(G,1)$
with finite $(n+1)$ skeleton $X^{{\prime}(n+1)}$ and let $X$
be the universal cover of $X^{\prime}$. The vector space of cellular
chains ${C_i}(X,{\mathbb{R}})$ is equipped with the $l_1$ norm
for a basis of $i$--cells. Then the boundary maps
${\delta_{i+1}} \co {C_{i+1}}(X,{\mathbb{R}}) \rightarrow {C_i}(X,{\mathbb{R}})$
are bounded linear and (owing to the finiteness of the $n+1$--skeleton)
one gets quasi-isometry invariant homology groups $H_i^{(1)}(X,{\mathbb{R}})$
for $i \leq n$. Since these homology groups are quasi-isometry invariant
it makes sense to define 
$H_i^{(1)}(G,{\mathbb{R}}) = H_i^{(1)}(X,{\mathbb{R}})$ for $i \leq n$
for any such $X$. The following Theorem of Gersten's occurs in \cite{gersten1}.

\begin{theorem}
The finitely presented group $G$ is hyperbolic if and only if
$H_1^{(1)} (G,{\mathbb{R}}) = 0$. Moreover, if $H$ is a finitely
generated subgroup of $G$ then $H$ is quasiconvex if and only if the map
$H_1^{(1)} (H,{\mathbb{R}}) \rightarrow H_1^{(1)} (G,{\mathbb{R}})$
induced by inclusion is injective.
\end{theorem}

Earlier results along these lines had been found in \cite{gersten},
\cite{gersten2}, \cite{gersten3}.

In a different direction, one would like a purely group-theoretic
characterization of quasiconvexity. We start with some definitions.

\begin{defn}
Let $H$ be a subgroup of a  group $G$. We say that
the elements $\{g_i |1 \le i \le n\}$ of $G$ are 
essentially distinct if $Hg_i \neq Hg_j$ for $i \neq j$.
Conjugates of $H$ by essentially distinct elements are called
 essentially distinct conjugates.
\end{defn}
Note that we are abusing notation slightly here, as a conjugate of $H$
by an element belonging to  the normalizer of $H$ but not belonging to
 $H$ is still essentially distinct from $H$.
Thus in this context a conjugate of $H$ records (implicitly) the conjugating
element.
\begin{defn} 
We say that the height of an infinite subgroup $H$ in $G$ is $n$  if  
there exists a collection of $n$ essentially distinct conjugates
of $H$ such that the intersection of all the elements of the collection is 
infinite  and $n$ is maximal possible. We define the height of a finite 
subgroup  to be $0$.
\end{defn}
The  main  theorem of \cite{GMRS} states:
\begin{theorem}
If $H$ is a quasiconvex subgroup of a hyperbolic group $G$,then
$H$ has finite height.
\label{gmrs}
\end{theorem}
The following question of Swarup was prompted partly by this result:
\par\medskip
{\bf Question}\stdspace (Swarup)\stdspace Suppose $H$ is a finitely presented
subgroup of a hyperbolic group $G$. If $H$ has finite height is $H$ quasiconvex in $G$? 
\par\medskip
So far only some partial answers have been obtained. The first result
is due to Scott and Swarup:

\begin{theorem}{\rm\cite{scottswar}}\qua 
Let $1 \rightarrow H \rightarrow G \rightarrow  {\mathbb{Z}} \rightarrow 1$
be an exact sequence of hyperbolic groups
induced by a pseudo Anosov  diffeomorphism of a
closed surface with
fundamental group $H$.
  Let $H_1$ be a finitely generated subgroup of infinite index in 
$H$. Then $H_1$ is quasiconvex in $G$.
\end{theorem}
In \cite{mitra4} an analogous result for free groups was derived. The methods
also provide a different proof of Scott and Swarup's theorem
above:
\begin{theorem}{\rm\cite{mitra4}}\qua  
 Let $1 \rightarrow H \rightarrow G \rightarrow \mathbb{Z} 
\rightarrow 1$ be an exact sequence of hyperbolic groups
induced by a hyperbolic automorphism $\phi$ of the free group
$H$. 
Let  $H_1 ( \subset H)$ be  a  finitely
generated distorted 
subgroup of $G$.
Then there exist  $N > 0$ and a free factor $K$ of $H$ such that
the conjugacy class of $K$ is preserved by $\phi^N$ and $H_1$
contains a finite index subgroup of a conjugate of $K$.
\end{theorem}
Another special case where one has a positive answer is
the following:
\begin{theorem}{\rm\cite{mitra5}}\qua 
 Let $G$ be a hyperbolic group splitting over $H$ (ie
$G = {G_1}{*_H}{G_2}$ or $G = {G_1}{*_H}$) with hyperbolic
vertex and edge groups. Further, assume
the  two inclusions of $H$ are  quasi-isometric embeddings.
 Then $H$ is of finite height in $G$
if and only if it is quasiconvex in $G$.
\label{mit5}
\end{theorem}
Swarup's question is therefore still open in the following 
special case, which can be regarded as  a next step following
the Theorems of \cite{mitra4} and \cite{mitra5} above. 
\par\medskip
{\bf Question}\stdspace Suppose  $G$ splits over $H$ satisfying the
hypothesis of Theorem \ref{mit5} above and $H_1$
is a quasiconvex subgroup of $H$. If $H_1$ has finite height in $G$
is it quasiconvex in $G$? More generally, if $H_1$ is an edge
group in a hyperbolic graph of hyperbolic groups satisfying the
qi--embedded condition, is $H$ quasiconvex in $G$ if and only if it
has finite height in $G$?
\par\medskip
A closely related problem  can be formulated in more geometric terms:

{\bf Question}\stdspace Let $X_G$ be a finite 2 complex with fundamental
group $G$. Let $X_H$ be a cover of $X_G$ corresponding to the 
finitely presented subgroup $H$. Let $I(x)$ be the injectivity
radius of $X_H$ at $x$.

Does $I(x) \rightarrow \infty$ as $x \rightarrow \infty$ imply that
$H$ is quasi-isometrically embedded in $G$?
\par\medskip
A positive answer to this question for $G$ hyperbolic would provide
a positive answer to Swarup's question.

The answer to this question is negative if one allows $G$
to be only  finitely generated instead of finitely presented
as the following example shows:
\par\medskip
{\bf Example}\stdspace Let $F = \{a,b,c,d \}$ 
denote the free group
on four generators.
Let ${u_i} = {a^i}{b^i} $ and ${v_i} = c^{f(i)}d^{f(i)}$
for some function $f \co \mathbb{N} \rightarrow \mathbb{N}$.
Introducing a stable letter $t$ conjugating $u_i$ to $v_i$
one has a finitely generated HNN extension $G$. The
free subgroup generated by $a, b$ provides a negative answer
to the question above for suitable choice of $f$. 
In fact one only requires that $f$ grows faster than any
linear function. 

If $f$ is
recursive one can embed the resultant $G$ in a finitely 
presented group by Higman's Embedding Theorem. But then 
one might lose malnormality of the free subgroup generated by $a, b$.
If one can have some control over the embedding in a finitely
presented group, one might look for a counterexample.
A closely related example was shown to the author by Steve Gersten.
\par\medskip

So far  the following question (attributed to Bestvina and Brady)
remains open:
\par\medskip
{\bf Question}\stdspace Let $G$ be a finitely presented group with a finite
$K(G,1)$. Suppose moreover that $G$ does not contain any subgroup
isomorphic to $BS(m,n)$. Is $G$ hyperbolic?
\par\medskip
A malnormal counterexample to Swarup's question would provide a 
counterexample for the above question (observed independently
by M. Sageev).

\section{Boundary extensions}

The purpose of this section is to take an asymptotic rather than
a coarse point of view and expose some of the problems from this perspective.
Since virtually all the work in this area involves actions of hyperbolic
groups on hyperbolic metric spaces we restrict ourselves mostly to
this.

Roughly speaking, one 
would like to know what happens `at infinity'.  We put this in the more
general context of a 
hyperbolic group $H$ acting freely and properly discontinuously
by isometries
on a proper hyperbolic metric space $X$. Then there is a natural map
 $i \co \Gamma_H \rightarrow X$, sending  the vertex set of $\Gamma_H$ to
the orbit of a point under $H$, and connecting images of adjacent vertices
in $\Gamma_H$ by geodesics in $X$. Let $\widehat{X}$ denote the Gromov
compactification of $X$.
 
 The basic question discussed in this section is the following:
\par\medskip
{\bf Question}\stdspace Does the continuous proper map $ i\co \Gamma_H
\rightarrow X$ extend to a continuous map $\hat i \co \widehat{\Gamma_H} \rightarrow \widehat{X}$?
\par\medskip
A measure--theoretic version of this question was asked by Bonahon
in \cite{Bonahon}. A positive answer to the above would imply
a positive answer to Bonahon's question. 
Related questions in the context of Kleinian groups have been studied by
Cannon and Thurston \cite{CT}, 
Bonahon \cite{Bonahon1}, Floyd \cite{Floyd} and  Minsky \cite{Minsky}.

Much of the work around this problem was inspired by a seminal
(unpublished)
paper of Cannon and Thurston \cite{CT}. The main theorem of \cite{CT}
states:

\begin{theorem}{\rm\cite{CT}}\qua
Let $M$ be a closed hyperbolic $3$--manifold fibering over the circle with 
fiber $F$. Let $\widetilde F$ and $\widetilde M$ denote the universal
covers of $F$ and $M$ respectively. Then $\widetilde F$ and $\widetilde M$
are quasi-isometric to ${\mathbb{H}}^2$ and ${\mathbb{H}}^3$ respectively. 
Let
${{\mathbb{D}}^2}={\mathbb{H}}^2\cup{\mathbb{S}}^1_\infty$ and 
${{\mathbb{D}}^3}={\mathbb{H}}^3\cup{\mathbb{S}}^2_\infty$
denote the standard compactifications. 
Then  the usual inclusion of $\widetilde F$ into $\widetilde M$
extends to a continuous map from ${\mathbb{D}}^2$ to ${\mathbb{D}}^3$.
\label{CT}
\end{theorem}

The proof of the above theorem involved the construction of a local
`Sol-like' metric using  affine structures on surfaces
coming from  stable and unstable
foliations. Coupled with Thurston's hyperbolization of $3$--manifolds
fibering over the circle one has a very explicit description
of the boundary extension. 

Using these (local) methods Minsky \cite{Minsky} generalized this
theorem to the following:

\begin{theorem}{\rm\cite{Minsky}}\qua
 Let $\Gamma$ be a 
 Kleinian 
group isomorphic (as a group) to the fundamental group of a closed
surface, such that ${{\mathbb{H}}^3}/{\Gamma} = M$ has injectivity radius 
uniformly bounded below by some $\epsilon > 0$. Then there exists a continuous
map from the Gromov boundary of $\Gamma$ (regarded as an abstract group)
to the limit set of $\Gamma$ in ${\mathbb{S}}^2_{\infty}$.
\label{minsky}
\end{theorem}

Finally Klarreich \cite{klarreich}
 generalized the above theorem to the case of freely indecomposable
Kleinian groups. A different proof was given by the author
 \cite{mitra3} (see below).

\begin{theorem}[\cite{klarreich},\cite{mitra3}]
 Let $\Gamma$ be a freely indecomposable
 Kleinian 
group, such that ${{\mathbb{H}}^3}/{\Gamma} = M$ has injectivity radius 
uniformly bounded below by some $\epsilon > 0$. Then there exists a continuous
map from the Gromov boundary of $\Gamma$ (regarded as an abstract group)
to the limit set of $\Gamma$ in ${\mathbb{S}}^2_{\infty}$.
\label{locconn}
\end{theorem}

 Klarreich proved  Theorem \ref{locconn} by combining her
Theorem \ref{klarreich} below  with Theorem \ref{minsky} above.

\begin{theorem}{\rm\cite{klarreich}}\qua
Let $X$ and $Y$ be proper, geodesic Gromov--hyperbolic spaces, ${H_\alpha}$
a collection of closed, disjoint path-connected subsets of $X$,
and $h\co X\rightarrow Y$ a quasi--Lipschitz map such that for every
$H_\alpha$, $h$ restricted to $H_\alpha$ extends continuously to 
 the boundary at infinity.  Suppose that
the following hold:

\begin{enumerate}\renewcommand{\labelenumi}{{\rm(\theenumi)}}

\item  The complement in $X$ of the sets $H_\alpha$ is open and
path-connected as also the complement of $h({H_{\alpha}})$ in $Y$.

\item There is some real number $k > 0$  such that the sets $H_\alpha$
 are all $k$--quasiconvex in $X$ and $h({H_{\alpha}})$'s are 
$k$--quasiconvex in $Y$.

\item There is a real number $c > 0$ such that $d({H_\alpha},{H_\beta}) > c$
and such that $d(h({H_\alpha}),h({H_\beta})) > 0$ for all $\alpha$ and $\beta$.
\end{enumerate}

Then if the map $h$ induced on the electric spaces is a quasi-isometry,
$h$ extends continuously to a continuous map from the boundary of
$X$ to the boundary of $Y$. Here the electric spaces are the spaces
obtained from $X$ and $Y$ by collapsing each space $H_\alpha$ (or
$h({H_\alpha})$) to points: they inherit path metrics from $X$ and $Y$.

\label{klarreich}
\end{theorem}

One should note  that since Cannon and Thurston's Theorem
\ref{CT}
deals with asymptotic behavior it might well be regarded as a theorem
in coarse geometry. The above Theorems are all of this form. But
the proof techniques in \cite{CT}, \cite{Minsky}
are local as they rely on Thurston's theory
of singular foliations of surfaces. In \cite{mitra1} and \cite{mitra3}
a different approach was described using purely large-scale
techniques giving generalized versions of Theorems \ref{CT}
\ref{minsky} and \ref{locconn}.

\begin{theorem}{\rm\cite{mitra1}}\qua
 Let $G$ be a hyperbolic group and let $H$ be a hyperbolic subgroup
that is normal in $G$. Let 
$i \co \Gamma_H\rightarrow\Gamma_G$ be the continuous proper
 embedding of $\Gamma_H$ in $\Gamma_G$ described above.
 Then $i$ extends to a continuous
map $\hat{i}$ from
$\widehat{\Gamma_H}$ to $\widehat{\Gamma_G}$.
\end{theorem}
A more useful generalization of Theorem \ref{CT} is:
\begin{theorem}{\rm\cite{mitra3}}\qua
 Let (X,d) be a tree (T) of hyperbolic metric spaces satisfying the
quasi-isometrically embedded condition.  Let $v$ be a vertex of $T$. Let
$({X_v},d_v)$ denote the hyperbolic metric space corresponding to $v$.
 If $X$ is hyperbolic then the inclusion of $X_v$ in $X$ extends
continuously to the boundary.
\label{tree}
\end{theorem}
A direct consequence of Theorem \ref{tree} above is the following:
\begin{cor}
 Let $G$ be a hyperbolic group acting cocompactly on a simplicial tree
$T$ such that all vertex and edge stabilizers are hyperbolic. Also
suppose that every inclusion of an edge stabilizer in a vertex stabilizer
is a quasi-isometric embedding. Let $H$ be the stabilizer of a vertex or
edge of $T$. Then an inclusion of the Cayley graph of $H$ into that of
$G$ extends continuously to the boundary.
\end{cor}
In \cite{BF}, Bestvina and Feighn give sufficient conditions for
a graph of hyperbolic groups to be hyperbolic. Vertex and edge
subgroups are thus natural examples of hyperbolic subgroups of hyperbolic
groups. 
These examples are covered by the above corollary.

Using Thurston's pleated surfaces technology one then gives a `coarse'
proof of Theorem \ref{locconn}. With some further work and using a theorem of Minsky \cite{Minsky2},
one can give \cite{mitra3} a `partly coarse' proof of another result of 
Minsky \cite{Minsky}:  Thurston's Ending Lamination Conjecture
for geometrically tame manifolds with 
freely indecomposable fundamental group and a uniform lower bound on
injectivity radius.

\begin{theorem}{\rm\cite{Minsky}}\qua Let $N_1$ and $N_2$ be homeomorphic
hyperbolic $3$--mani-folds with freely indecomposable fundamental
group. Suppose there exists a uniform lower bound $\epsilon > 0$  on
the injectivity radii of $N_1$ and $N_2$. If the end invariants of
corresponding ends of $N_1$ and $N_2$ are equal, then $N_1$ and $N_2$ are isometric.
\end{theorem}

One should note here that the coarse techniques referred to 
circumvent only the building of a `model manifold' --- a local
construction in \cite{Minsky}.
It might be worthwhile to obtain a coarse proof of the main
theorem of \cite{Minsky2}. A positive answer to the following
coarse question will do the job (as can be seen from \cite{mitra3}):
\par\medskip
{\bf Question}\stdspace Let $\sigma \co \mathbb{N} \rightarrow Teich(S)$ be a 
map. For $l$ a closed curve on $S$, let $l_i$ denote the length 
of the shortest curve freely homotopic to $l$ on ${\sigma}(i)$.

Suppose there exists $\lambda > 1$ such that for all closed
curves   $l$ on $S$ one has

\begin{center}

${\lambda}l_i \leq max({l_{i-1}},{l_{i+1}})$ for all $i \in
\mathbb{N}$.

\end{center}

Then does $\sigma$ lie in a bounded neighborhood of a Teichmuller
geodesic?
\par\medskip
The above question was motivated in part by the `hallways flare'
condition of \cite{BF} and a recent relative hyperbolicity
result of Masur--Minsky \cite{masur-minsky}. 

Since a continuous image of a compact locally connected set is locally
connected \cite{hock-young} Theorem \ref{locconn} also shows that the
limit sets of freely indecomposable Kleinian groups with a uniform
lower bound on the injectivity radius are locally connected. The issue
of local connectivity has received a lot of attention lately due to
some recent foundational work of Bowditch and Swarup \cite{bow1},
\cite{bow2}, \cite{bow3}, \cite{bow4}, \cite{swarup} following earlier
work by Bestvina and Mess \cite{bes-mess}.

\begin{theorem}[\cite{bow1}, \cite{swarup}]Let $H$ be a one-ended hyperbolic
group. Then its\break boundary is locally connected. Next assume $H$ does
not split over any two-ended group and acts on a proper hyperbolic
metric space $X$ with limit set $\Lambda \subset \partial X$. Then
$\Lambda$ is locally connected.
\label{bow}
\end{theorem}

The existence of continuous boundary extensions in general
would thus imply (using Theorem \ref{bow}) local connectivity
of limit sets of hyperbolic groups acting on proper
hyperbolic metric spaces. One wonders if some kind of a converse
exists. 

Such speculations are prompted on the one hand by
Theorem \ref{bow} and  by the following observation.
Let
$\Gamma$ be a simply degenerate Kleinian group isomorphic to
a surface group. Further assume $\Gamma$ has no parabolics.
Let $\Lambda$ be the limit set of $\Gamma$, $\Omega$
its domain of discontinuity and $X$ the boundary of the
convex hull of $\Lambda$. Then `nearest point projections' give 
a natural homeomorphism between $\Omega$ and $X$. From this it is
easy to conclude that a continuous boundary extension exists if
and only if a neighborhood of $\Lambda$ in $S^2_\infty$
deformation retracts
onto $\Lambda$. In this special case therefore local connectivity
is equivalent to continuous boundary extensions. 

Before concluding this section it is worth pointing out that one needs
finer invariants than distortion to understand asymptotic extrinsic
geometry. One way of approaching the problem is to consider extrinsic
geometry of rays (starting at $1 \in \Gamma_H$) and describe those
which are not quasigeodesics in the ambient space $X$.  If one looks at
bi-infinite geodesics instead of rays one gets `ending laminations'.
For $3$--manifolds fibering over the circle with fiber $F$ and monodromy
$\phi$ one can think of these as the stable and unstable foliations
of $\phi$. Motivated by this, the author gave a more group theoretic
description in \cite{mitra2} in the special case of a hyperbolic normal
subgroup of a hyperbolic group.

Recall that for a hyperbolic $3$--manifold $M$ fibering over the circle with
fiber $F$  Cannon and Thurston
show in \cite{CT} 
that the usual inclusion of $\widetilde F$ into $\widetilde M$
extends to a continuous map from ${\mathbb{D}}^2$ to ${\mathbb{D}}^3$. An
explicit description of this map was also described in \cite{CT} in
terms of `ending laminations' [See \cite{Thurstonnotes} for definitions].
The explicit description depends on Thurston's theory of stable and unstable
laminations for pseudo-anosov diffeomorphisms of surfaces \cite{FLP}.
In the case of  normal hyperbolic subgroups of hyperbolic groups, 
though existence of a continuous extension 
${\hat{i}} \co \widehat{\Gamma_H} \rightarrow \widehat{\Gamma_G}$
was proven in \cite{mitra1}, an explicit description was missing.
In \cite{mitra2}
 some parts  of Thurston's theory of ending laminations
were generalized  to the context of
normal hyperbolic subgroups of hyperbolic groups. Using this 
an explicit  description of the continuous boundary
extension 
$\hat{i} \co \widehat{\Gamma_H} \rightarrow \widehat{\Gamma_G}$
was given for $H$
a normal hyperbolic 
subgroup of a  hyperbolic group $G$.

In general,
if

\begin{center}

$1 \rightarrow H \rightarrow G \rightarrow Q \rightarrow 1$

\end{center}
is an exact sequence of finitely presented groups where $H$, $G$
and hence $Q$ (from \cite{Mosher}) are hyperbolic, one has   
ending laminations 
naturally parametrized by points in the boundary $\partial\Gamma_Q$
of the quotient
group $Q$.

Corresponding to every element $g\in G$ there exists an automorphism 
of $H$ taking $h$ to $g^{-1}hg$ for $h\in H$. Such an automorphism induces
a bijection $\phi_g$ of the vertices of $\Gamma_H$. This gives rise to a map
from $\Gamma_H$ to itself, sending an edge [$a,b$] linearly to a shortest
edge-path joining $\phi_g (a)$ to $\phi_g (b)$. 

Fixing $z\in{\partial{\Gamma_Q}}$ for the time being
(for notational convenience) we shall define the
set of ending laminations corresponding to $z$.

Let $[1,z)$ be a semi-infinite geodesic ray in $\Gamma_Q$ starting at the
identity $1$ and converging to $z\in{\partial}{\Gamma_Q}$. Let
$\sigma$ be a single-valued quasi-isometric section of $Q$ into $G$.
Let $z_n$ be the vertex on $[1,z)$ such that ${d_Q}(1,{z_n}) = n$ and let
${g_n} = {\sigma}({z_n})$. 

Given $h\in{H}$ let ${\Sigma}_n^h$ be the ($H$--invariant) collection of
all free homotopy representatives (or shortest representatives in the
same conjugacy class) of  ${\phi}_{g_n^{-1}}(h)$ in $\Gamma_H$.  
Identifying equivalent geodesics in ${\Sigma}_n^h$ one obtains a subset
$S_n^h$ of (unordered) pairs of points in 
${\widehat{{\Gamma}_H}}$. The intersection 
with ${\partial}^{2}{\Gamma_H}$ of the union
of all subsequential limits (in the Chabauty topology)
of $\{{S_n^h}\}$ will be denoted by 
${\Lambda}_z^h$. 

\begin{defn} The set of {\it ending laminations corresponding to
$z\in{\partial}{\Gamma_Q}$ } is given by
$$\Lambda_z = {{\bigcup}_{h\in{H}}}{{\Lambda}_z^h}.$$
\end{defn} 

\begin{defn}The set $\Lambda$ 
of all {\it ending laminations} is defined
by 
$$\Lambda = {{\bigcup}_{z\in{\partial}{\Gamma_Q}}}{{\Lambda}_z}.$$
\end{defn}

It was shown in \cite{mitra2}  that the continuous boundary extension
$\hat{i}$ identifies end-points of 
leaves of the ending lamination. Further 
 if $\hat{i}$ identifies
a pair of points in $\partial\Gamma_H$, then a bi-infinite geodesic
having these points as its end-points is a leaf of the ending lamination.

Similar descriptions of laminations 
 have been used by Bestvina, Feighn and Handel for free groups
\cite{BFH}. Using these two descriptions in conjunction gives
further information eg about subgroup structure  \cite{mitra4}.

\section{Other invariants in extrinsic geometry}
To fix notions consider a finitely generated group $H$ acting on a
path-metric space $X$. As mentioned in the introduction distortion
arises out of comparing the intrinsic metric on $\Gamma_H$ to
the metric inherited from the ambient space $X$. Alternately this
can be regarded as arising out of comparing filling functions,
where one fills a copy of $S^0$ in $\Gamma_H$ and $X$ and compares
the sizes of the chains required.

In Chapter 5 of \cite{Gromov2} Gromov defines several filling 
invariants of spaces. Each of these gives rise to a relative
version and corresponding distortion functions. Recall some
of these from \cite{Gromov2}.

Given a simplicial $n$--cycle $S$ in a homotopically (or homologically)
$n$--connected simplicial complex $X$ one constructs fillings of $S$
by $(n+1)$ chains in $X$.
\begin{defn}
Filling volume, denoted $FillVol_n (S,X)$ is the infimal simplicial
volume of $(n+1)$ chains filling $S$.
\end{defn}
\begin{defn}
Filling radius, denoted $FillRad_n (S,X)$ is the minimal
$R$ such that $S$ bounds in an $R$--neighborhood ${U_R}(S) \subset X$.
\end{defn}

A host of other filling invariants are defined in \cite{Gromov2}
but we focus on these two. 

We will define relative versions of the above two notions. Since
the definitions of these invariants require $n$--connectedness
of the spaces we shall assume that whenever these invariants are
defined, the spaces in question are quasi-isometric to (or admit
thickenings that are) $n$--connected. It will be clear that one gets
quasi-isometry invariants in the process. Reference to an explicit
quasi-isometry may at times be suppressed.
  
Distortion of 
$Fillvol_n$ and $FillRad_n$ can be defined in a somewhat more general
context. Fix  classes ${{\mathcal{S}}_n}(X)$ and
${{\mathcal{S}}_n}(Y)$ 
of $n$--cycles 
in $X$, $Y$ respectively (eg one might restrict to connected cycles
or images of spheres) such that
${{\mathcal{S}}_n}(X) \subset {{\mathcal{S}}_n}(Y)$. Let $f_n$ be one
be one of the functions $FillVol_n$ or $FillRad_n$. Define

\begin{center}
${{\mathcal{S}}_n}({f_n},m,X) $ =
$\{ S \in {\mathcal{S}}_n (X) : {f_n}(S) \leq m \}$.
\end{center}

Finally define

\begin{center}

$Disto(f,X,Y,m) = sup ({f_n}(S,Y))$ \\
where the $sup$ is taken over 
$S \in {{\mathcal{S}}_n}({f_n},m,X) \cap {{\mathcal{S}}_n}(Y)$.

\end{center}

For $n = 0$, ${\mathcal{S}} -$ the set of maps of the $0$--sphere $S^0$
and $f_0 = FillVol_0$ or $FillRad_0$ we get back the original
distortion function. Note that $FillRad_0$ is approximately half
of $FillVol_0$. 

For $n = 1$, ${\mathcal{S}} $ the set of maps of  $S^1$
and $f_0 = FillVol_1$ we get {\it area distortion} in the sense
of Gersten \cite{gerarea}. Distortion has been surveyed in Section 2.
We give a brief sketch of Gersten's results on area distortion.

\begin{defn} An automorphism of a finitely presented group
is {\it tame} if it lifts to an automorphism of the free group
on its generators, preserving the normal subgroup generated by relators.
\end{defn}
\begin{theorem}
Let $\phi$ be a tame automorphism of a one-relator group $G$. Then
area is undistorted for $G \subset G \rtimes_{\phi}\mathbb{Z}$.
\end{theorem}
In \cite{gerarea} Gersten shows that in
extensions of $\mathbb{Z}$ by finitely presented groups $G$
area distortion of $G$ is at most an exponential of an
isoperimetric function for the extension. Moreover, he describes
examples of undistorted (in the usual sense of length) subgroups
that exhibit area distortion. He observes further that for torus
bundles over the circle with
Sol geometry, area in  the fiber subgroup 
is undistorted whereas length is exponentially distorted.

Gersten showed further that area is undistorted for finitely
presented subgroups of finitely presented groups of cohomological
dimension 2. From this it follows that finitely presented
subgroups $H$ of hyperbolic groups $G$ are finitely presented
provided $G$ has cohomological dimension 2  or $G$ is a hyperbolic
small cancellation group \cite{gersubgp}.

The remaining distortion functions are yet to be studied
systematically.
The first class of examples where $Disto({FillVol_n},X,Y,m)$
seem interesting and tractable are examples coming from extensions of 
$\mathbb{Z}$ by ${\mathbb{Z}}^n$, ie for
$G = {\mathbb{Z}}^n {\rtimes}_{\phi}{\mathbb{Z}}$ where 
$\phi \in GL_n{\mathbb{Z}}$. Such examples have been
studied by  Bridson \cite{bridson1}
and Bridson and Gersten \cite{brid-ger}.

Much less is known about $Disto({FillRad_n},X,Y,m)$. These functions
are related to topology of balls in groups (Chapter 4 of
\cite{Gromov2}).
For a group $\Gamma$ admitting a uniformly $k$--connected thickening
$X$ (see \cite{Gromov2} for definitions) Gromov defines 
${\overline{R_k}}(r)$ 
to be the infimal radius $R \geq r$
such that the inclusion of balls $B(r) \subset B(R)$ is 
$k$--connected.

The following observations are straightforward generalizations of 
corresponding statements (for $n = 0$) on pages 74--76 of
\cite{Gromov2}. Fix a group $\Gamma^{\prime}$ and a subgroup $\Gamma$.

\begin{prop}
If $Disto({FillRad_n}, \Gamma^{\prime} , \Gamma , m)$ is
superexponential
in $m$ then the function ${\overline{R_k}}(m)$  for
$({\Gamma},dist_{\Gamma^{\prime}}|{\Gamma})$
grows faster than any linear function $Cm$.
\end{prop}

\begin{prop}
Take two copies of $({\Gamma^{\prime}},{\Gamma \subset
\Gamma^{\prime}})$
and let ${\Gamma_1} = {\Gamma^{\prime}}*_{\Gamma}{\Gamma^{\prime}}$
be the double. Then the function 
${\overline{R_k}}(m)$ 
 for
${\Gamma}_1$ is minorized by 
${\overline{R_{k-1}}}(m)$ 
for 
$({\Gamma},dist_{\Gamma^{\prime}}|{\Gamma})$.
\label{ri}
\end{prop}

This leads to the following

\par\medskip
{\bf Question}\stdspace Do there exist pairs of groups $H \subset G$
(with $n$--connected inclusions of thickenings of the Cayley 
Graph) such that  $Disto({FillRad_n}, \Gamma_G , \Gamma_H , m)$ is
superexponential
in $m$?
\par\medskip

A positive answer will furnish 
(via Proposition \ref{ri}) examples of groups with fast growing
${\overline{R_k}}(m)$ 
for $k \geq 2$ (page 80 of \cite{Gromov2}).
No such example has been found yet.

\Addresses\recd

\end{document}